\documentclass{amsart}
\usepackage{latexsym}
\usepackage{amssymb}
\usepackage{amsmath}
\usepackage{amsfonts}

\def\R{\mathbb{R}}

\newcommand{\bmat}{\left[\begin{matrix}}
\newcommand{\emat}{\end{matrix}\right]}
\usepackage[latin1]{inputenc} 
\usepackage{amsthm}
\usepackage{amsxtra}
\usepackage{amscd}
\usepackage{setspace}   
\usepackage{mathrsfs}   

\newtheorem{theorem}{Theorem}

\theoremstyle{remark}

\theoremstyle{definition}

\newcommand{\Z}{\mathbb{Z}}

\usepackage{xcolor}
\usepackage{enumerate}

    \title{A new upper bound on the smallest counterexample to the Mertens conjecture}
    \author{John Rozmarynowycz and Seungki Kim}

\begin{document}
\maketitle

\begin{abstract}

We report the finding of the new upper bound on the lowest positive integer $x$ for which the Mertens conjecture
\begin{equation*}
\left| \sum_{1 \leq n \leq x} \mu(n) \right| < \sqrt{x}
\end{equation*}
fails to hold: $x < \exp(1.017 \times 10^{29})$, an improvement over previously known $\exp(1.59 \times 10^{40})$ due to Kotnik and te Riele \cite{KtR}.

\end{abstract}


\section{Introduction}

Perhaps one of the most striking application of the LLL reduction algorithm (\cite{LLL}; also see \cite{LLLbook}) to number theory is the disproof of the Mertens conjecture by Odlyzko and te Riele \cite{OtR} in 1985. It is a conjecture made by Mertens \cite{Mertens} in 1897 stating that
\begin{equation*}
M(x) := \left| \sum_{1 \leq n \leq x} \mu(n) \right| < \sqrt{x} \mbox{ for any $x > 1$,}
\end{equation*}
where $\mu(n)$ is the usual M\"obius function.

This conjecture lasted for nearly a century, until in \cite{OtR} the authors solved the associated problem in simultaneous diophantine approximation, by reformulating it as a lattice reduction problem and applying the then-cutting-edge LLL algorithm. There are two natural follow-up questions to ask, both still open today:
\begin{enumerate}[(i)]
\item What is the correct asymptotic growth rate of $M(x)$?
\item What is the smallest $x$ for which $|M(x)| \geq \sqrt{x}$?
\end{enumerate}
The present paper focuses on the latter question. A theorem of Pintz \cite{Pintz} provides an approach via diophantine approximation again, based on which Kotnik and te Riele \cite{KtR} showed in 2006
\begin{equation*}
x < \exp(1.59 \times 10^{40}),
\end{equation*}
which remains the best known bound to this date. A numerical study of Kotnik and van de Lune \cite{KvdL} conjectures that the smallest counterexample should be of size about $\exp(5.15 \times 10^{23})$.

Meanwhile, there have been huge and rapid improvements in the art of lattice reduction, largely motivated by post-quantum cryptography. LLL itself has seen several improvements which led to substantial speedups in practice (e.g. \cite{NS05}; see also \cite{RH}). Furthermore, much stronger algorithms have been proposed. Especially noteworthy is the BKZ algorithm, originally due to Schnorr and Euchner \cite{SE94}, that has undergone a series of optimizations in both output quality and complexity since the last decade (e.g. \cite{CN10}, \cite{AWHT16}). Nowadays, fpLLL's implementation \cite{fpLLL} of BKZ yields a result of strength --- e.g. the quality of the diophantine approximation --- unachievable by LLL in the blink of an eye on a personal laptop.

Our idea was simply to adopt the more recent and powerful lattice reduction in place of LLL. By randomized runs of BKZ on a personal laptop, guided by some common-sense knowledge on lattices, we obtained the bound
\begin{equation*}
x < \exp(1.017 \times 10^{29})
\end{equation*}
on the smallest counterexample to the Mertens conjecture, bringing us closer to the conjectured $x \approx \exp(5.15 \times 10^{23})$. Furthermore, we found a number of suggestive data points in this conjectured range --- see Section 3 below.

This paper demonstrates only a tiny portion of the current art of lattice reduction. In a forthcoming work, we employ much more powerful tools and techniques for the goal of attaining the conjectured bound of \cite{KvdL} and perhaps even further. We hope our work motivates more applications of the recent advances in the computational lattice problems to number theory.

\subsection*{Acknowledgments}

Both authors acknowledge support from the NSF grant CNS-2034176. We thank Yuntao Wang for his help with the powerful progressive BKZ, and Phong Nguyen for helpful discussions.

\section{Outline of the approach}

We denote by $\mu(n)$ and $\zeta(s)$ the usual M\"obius function and the Riemann zeta function, respectively. $\rho$ denotes a zero of $\zeta(s)$ with $\mathrm{Re}\,\rho = \frac{1}{2}$, and for a given $\rho$ we denote $\gamma := \mathrm{Im}\,\rho, \alpha := |\rho\zeta'(\rho)|^{-1}, \psi := \mathrm{arg}(\rho\zeta'(\rho))$. There are two different ways in the literature to index the zeroes and the associated quantities: $\{\rho_i\}, \{\gamma_i\}, \{\alpha_i\}, \{\psi_i\}$ are ordered so that $\gamma_i < \gamma_{i+1}$ for all $i$, and $\{\rho^*_i\}, \{\gamma^*_i\}, \{\alpha^*_i\}, \{\psi^*_i\}$ are ordered so that $\alpha^*_i > \alpha^*_{i+1}$ for all $i$.

Our approach, as with \cite{KtR}, is based on the following result.
\begin{theorem}[Pintz \cite{Pintz}] \label{thm:pintz}
Let
\begin{equation} \label{eq:hp}
h_P(y) := 2\sum_{\gamma < 14000} \alpha\exp(-1.5 \cdot 10^{-6}\gamma^2)\cos(\gamma y - \psi).
\end{equation}
If there exists $y \in [e^7, e^{50000}]$ with $|h_P(y)| > 1+e^{-40}$, then $M(x) > \sqrt{x}$ for some $x < \exp(y + \sqrt{y})$.
\end{theorem}

The idea, due to \cite{OtR}, is that, to make $h_P$ as large as possible, one tries to minimize those $\gamma y - \psi$ (mod $2\pi$) with large weights, since \eqref{eq:hp} is approximately
\begin{equation} \label{eq:approx_hp}
2\sum_{i=1}^n \alpha^*_i\left(1-\left(\gamma^*_i y - \psi^*_i \mbox{ mod $2\pi$}\right)^2\right) + \mbox{(``error'')}
\end{equation}
for some $n \leq 100$, say. This leads them to consider the lattice in $\R^{n+2}$ generated by the columns of the matrix
\begin{equation} \label{eq:thematrix}
\begin{pmatrix}
-\lfloor \sqrt{\alpha^*_1}\psi^*_12^\nu \rfloor & \lfloor \sqrt{\alpha^*_1}\gamma^*_12^{\nu-10} \rfloor & \lfloor 2\pi\sqrt{\alpha^*_1}2^\nu \rfloor & 0  & \cdots & 0 \\
-\lfloor \sqrt{\alpha^*_2}\psi^*_22^\nu \rfloor & \lfloor \sqrt{\alpha^*_2}\gamma^*_22^{\nu-10} \rfloor & 0 & \lfloor 2\pi\sqrt{\alpha^*_2}2^\nu \rfloor & \cdots & 0 \\
\vdots & \vdots & \vdots & \vdots & \ddots & \vdots \\
-\lfloor \sqrt{\alpha^*_n}\psi^*_n2^\nu \rfloor & \lfloor \sqrt{\alpha^*_n}\gamma^*_n2^{\nu-10} \rfloor & 0 & 0 & \cdots & \lfloor 2\pi\sqrt{\alpha^*_n}2^\nu \rfloor \\
2^\nu n^4 & 0 & 0 & 0 & \cdots & 0 \\
0 & 1 & 0 & 0 & \cdots & 0
\end{pmatrix},
\end{equation}
where $n$ is as earlier, and $\nu$ is the parameter controlling the round-off of the involved quantities. In a reduced (by LLL or other algorithms) basis of \eqref{eq:thematrix}, one finds a vector, for some integers $p_i$'s and $z$,
\begin{align*}
&(p_1\lfloor 2\pi\sqrt{\alpha^*_1}2^\nu \rfloor + z \lfloor \sqrt{\alpha^*_1}\gamma^*_12^{\nu-10} \rfloor -\lfloor \sqrt{\alpha^*_1}\psi^*_12^\nu \rfloor, \ldots \\
&\ldots, p_n\lfloor 2\pi\sqrt{\alpha^*_n}2^\nu \rfloor + z \lfloor \sqrt{\alpha^*_n}\gamma^*_n2^{\nu-10} \rfloor -\lfloor \sqrt{\alpha^*_n}\psi^*_n2^\nu \rfloor, \pm2^\nu n^4, z)^\mathrm{tr},
\end{align*}
since $2^\nu n^4$ was chosen to be much larger than the rest of the entries of \eqref{eq:thematrix}. We let $y =  \pm z2^{-10}$, the sign being that of the second last entry; then the first $n$ entries can be seen to provide the minimizations of $\gamma^*_i y - \psi^*_i$ (mod $2\pi$) weighted by $\sqrt{\alpha_i^*}$.

One can also aim for a large negative value of $h_P$ by a similar construction, in which we replace $\psi^*$ in \eqref{eq:thematrix} by $\psi^* + \pi$.

From a slightly different perspective, what the reduction of \eqref{eq:thematrix} achieves is a solution to a certain \emph{approximate closest vector problem} (aCVP), that is, finding a vector of the lattice in $\R^{n+1}$ generated by the columns of 
\begin{equation} \label{eq:thematrix'}
\begin{pmatrix}
\lfloor \sqrt{\alpha^*_1}\gamma^*_12^{\nu-10} \rfloor & \lfloor 2\pi\sqrt{\alpha^*_1}2^\nu \rfloor & 0  & \cdots & 0 \\
\lfloor \sqrt{\alpha^*_2}\gamma^*_22^{\nu-10} \rfloor & 0 & \lfloor 2\pi\sqrt{\alpha^*_2}2^\nu \rfloor & \cdots & 0 \\
\vdots & \vdots & \vdots & \ddots & \vdots \\
\lfloor \sqrt{\alpha^*_n}\gamma^*_n2^{\nu-10} \rfloor & 0 & 0 & \cdots & \lfloor 2\pi\sqrt{\alpha^*_n}2^\nu \rfloor \\
1 & 0 & 0 & \cdots & 0
\end{pmatrix}
\end{equation}
that is very close to the ``target vector''
\begin{equation*}
\mathbf t := (-\lfloor \sqrt{\alpha^*_1}\psi^*_12^\nu \rfloor, \ldots, -\lfloor \sqrt{\alpha^*_n}\psi^*_n2^\nu \rfloor, 0)^\mathrm{tr} \in \R^{n+1}.
\end{equation*}
The matrix \eqref{eq:thematrix} is designed so that reducing it will execute what is nowadays known as \emph{Babai's nearest plane algorithm} \cite{Babai} (see also \cite[Chapter 2]{Prest}), one of the main approaches to aCVP to this day. The output quality, i.e. the distance from the found lattice vector to $\mathbf t$, is affected by the strength of the reduction algorithm used.

Let
\begin{equation*}
D = \prod_{i=1}^n \lfloor 2\pi\sqrt{\alpha^*_i}2^\nu \rfloor,
\end{equation*}
the determinant of \eqref{eq:thematrix'}. We expect a cube in $\R^{n+1}$ of volume $D$ to contain one lattice vector on average. Hence one expects a lattice vector $\mathbf v$ whose $L^{\infty}$-norm distance from $\mathbf t$ is at most about $\frac{1}{2}D^{\frac{1}{n+1}}$. Since the mass of a cube in a large dimension is concentrated on its corners, we expect each coordinate of $\mathbf v - \mathbf t$ to be of size $\approx \frac{1}{2}D^{\frac{1}{n+1}}$. Provided our reduction algorithm is strong enough, such $\mathbf v$ can indeed be found. In our experiments, we observed this heuristic to be correct up to a factor of a few hundreds, while $D^{\frac{1}{n+1}}$ is of magnitude $10^{30}$ to $10^{40}$.

This observation leads to an estimate on the expected size of the intended main term in \eqref{eq:approx_hp}: each term is of size $O_n(2^{-\frac{2\nu}{n+1}})$. This partially explains the first inequality in the condition $2n \leq \nu \leq 4n$ that \cite{OtR} imposed in their experiments. (The second inequality seems unnecessary to us though.) It also gives a heuristic estimate on the size of $y$ as a function of $\nu$; a convenient simple rule we found to work well in practice is $\log_{10} y \sim \log_{10}\nu - 5$ or $-6$.

\section{Experiment and result}

In search of the smallest $y$ for which $h_P(y) > 1 + e^{-40}$, we reduced lattices of the form \eqref{eq:thematrix} for $105 \leq \nu \leq 125$, $2n \leq \nu \leq 4n$, using BKZ-$\beta$ for $\beta = 20, 30, 40$ --- briefly speaking, higher $\beta$ makes the algorithm stronger and costlier. For each of these parameter choices, we took $500$ randomized bases of \eqref{eq:thematrix} to reduce, and for the randomization, we used the method of \cite{NS06} (see also \cite{BM} for more on basis randomization). For each $\nu$ and $\beta$, the computations took about one to two days on a personal laptop. We also ran the same experiment looking for a large negative value of $h_P$. 

Randomizing the input also randomizes the output to some extent, and the hope is that at least one of them yields $|h_P(y)| > 1 + e^{-40}$. In fact, this turned out to be a far more efficient strategy than applying a single high-quality (and high-cost) reduction to \eqref{eq:thematrix}. The reason is that a near-optimal solution to the associated aCVP problem, while taking a disproportionately longer time to compute, does not necessarily lead to a higher value of $h_P$. While the sum in \eqref{eq:approx_hp} is controllable by lattice reduction, the ``error'' part is not, and yet its size may fluctuate large enough to affect the outcome, either in or against our favor.

The values $\gamma, \alpha, \psi$ are taken from the data made public by Hurst \cite{Hurst} who computed them up to almost $10000$ decimal digits of precision, for which we are grateful. Our computations were made with 1024 binary digits of precision; it is easy to show that this is more than enough to estimate \eqref{eq:hp} well enough for our purpose. For the values presented below, they were checked again with 16384 binary digits of precision. The source code is made available on the second-named author's website: https://sites.google.com/view/seungki/

The lowest working value of $y$ is found with $\nu = 112, n = 53, \beta = 20$, in which
\begin{equation*}
y = 1017256208\ 7569945816\ 8018857216.806640625, h_P(y) = 1.0034372\ldots
\end{equation*}
By Theorem \ref{thm:pintz}, this shows that the first counterexample $x$ to the Mertens conjecture satisfies $x < \exp(1.017 \times 10^{29})$.

The conjecture of \cite{KvdL} that $x \approx \exp(5.15 \times 10^{23})$ corresponds to $\nu \approx 95$. We also found several very suggestive data points in this range (all with $\beta = 40$):
\begin{align*}
&y = 7272\ 5861306259\ 2936179649.8388671875, h_P(y) = 9.6027706\ldots \mbox{ for $\nu = 95, n = 50$} \\
&y = 3276\ 1262680303\ 1941538273.2919921875, h_P(y) = 9.6084449\ldots \mbox{ for $\nu = 95, n = 56$} \\
&y = 258\ 4924462692\ 5200109819. 8173828125, h_P(y) = -9.5313433\ldots \mbox{f or $\nu = 95, n = 61$}\\ 
&y = 5714\ 9077379396\ 8420303581, h_P(y) = -9.6006767\ldots \mbox{ for $\nu = 95, n = 64$}\\
&y = -81\ 4638194152\ 4511993798.2373046875, h_P(y) = -9.7588934\ldots \mbox{ for $\nu = 96, n = 58$}
\end{align*}

These examples encourage a further investigation; it is currently in progress.

\end{document}